\documentclass[12pt]{article}
\usepackage{amsmath}
\usepackage{amsfonts}

\title{On  the proof of some theorem on locally nilpotent subgroups in division rings }

\author{Bui Xuan Hai 
\\{\small\em Faculty of Mathematics $\&$ Computer Science, University of  Science}\\{\small\em VNU - HCM City}\\ {\small\em 227 Nguyen Van Cu str., Dist. 5, Ho Chi Minh City, Vietnam}\\ {\small\em e-mail: bxhai@hcmuns.edu.vn}\\
Nguyen Van Thin\\
{\small\em Faculty of Mathematics $\&$ Computer Science, University of Science}\\{\small\em VNU - HCM City}\\ {\small\em 227 Nguyen Van Cu str., Dist. 5, Ho Chi Minh City, Vietnam}\\ {\small\em e-mail: ngvthin@yahoo.com.vn}}

\date{}
\rm\em 
\normalsize
\begin{document}
\maketitle
\newcommand{\dpcm}{ \hfill \rule{3mm}{3mm}}
\def\Box{\dpcm}
\def\xd{\linebreak} 
\begin{abstract} 
 In  Hai-Thin (2009), there is a theorem,  stating that every  locally nilpotent subnormal subgroup in a division ring $D$  is central (see Hai-Thin (2009, Th. 2.2)).  Unfortunately, there is some mistake in the proof of this theorem. In this note we give  the another  proof of this theorem. 
\end{abstract}



\newpage

In  Hai-Thin (2009, Th. 2.2) there is the following theorem:\\

 \noindent
{\bf Theorem 1.} (see Hai-Thin (2009, Th. 2.2)) {\em Let $D$ be a division ring with the center $Z(D)$. Then, every locally nilpotent subnormal subgroup of $D^*$ is central, i.e. it is contained in $Z(D)$.}   \\   

Unfortunately,  in the proof of this theorem  there is some mistake. However, the theorem remains always true. Here  we give the another proof of this theorem. \\

\noindent
{\bf Lemma 2.} {\em Let $D$ be a non-commutative division ring and suppose that $G$ is a subnormal locally nilpotent  subgroup of $D^*$. If $x, y\in G$ with $[x, y]:=x^{-1}y^{-1}xy\neq 1$, then $[x, y]\not\in Z(D)$.}\\

\noindent
{\em Proof.} Since $G$ is subnormal in $D^*$, there exists the following series of subgroups:
$$G=G_n\triangleleft G_{n-1}\triangleleft\ldots\triangleleft G_1\triangleleft G_0=D^*.$$

Suppose that $x, y\in G$ with $[x, y]\neq 1$ and $x, y\in Z(D)$. Set
$$x_1:=[x+1, y], x_{i+1}:=[x_i, y], \forall i \geq 1.$$

We shall prove by induction that $x_k\in G_k, \forall k\in\{1, \ldots, n\}$. In fact, since $y\in G\leq G_k, \forall k$, it follows  that $x_1=(1+x)^{-1}y^{-1}(1+x)y\in G_1$. Now, suppose that $x_k\in G_k$, then, $x_{k+1}=[x_k, y]=(x_k^{-1}y^{-1}x_k)y\in G_{k+1}$. Thus, $x_k\in G_k, \forall k\in\{1, \ldots, n\}$. In particular, $x_n\in G_n=G.$  By supposition, $G$ is locally nilpotent, so the subgroup $\langle x_n, y\rangle$ generated by elements $x_n$ and $y$, is nilpotent. Hence, there exists some integer $t$ such that $x_t=1$. Let $P$ be a prime subfield of $D$. Then, by Scott(1987,  14.3.2, p. 432), there exists some polynomial $g\in P(b)[X]$ such that $g(x)=0$ and, moreover, $g$ is  not dependent either on $x$ or $y$, but it is dependent only on $b$. Clearly, for any $m\in \mathbb{Z}$ such that $mx\neq 0$, we have $[mx, y]=[x,y]=b$. Therefore, $g(mx)=0, \forall m\in \mathbb{Z}$ such that $mx\neq 0$. Evidently, this forces $Char(D)=p > 0$. Suppose that $f(X)=\displaystyle\sum_{i=0}^r a_iX^i\in P(b)[X]$ with $a_r=1$, is the minimal polynomial of the element $x$ over $P(b)$. Then
$$f(x)=\displaystyle\sum_{i=0}^r a_ix^i=0 \mbox{ and } y^{-1}f(x)y=\displaystyle\sum_{i=0}^r a_ib^ix^i=0.$$

It follows that $a_1(b-1)+ \ldots +a_r(b^r-1)x^{r-1}=0$. Therefore, $b^r-1=0$ or $b^r=1$. Hence $P(b)$ is a finite field. Thus, $x$ is algebraic over a finite field $P(b)$. By symmetry, we can conclude that $y$ is algebraic over $P(b)$ too. Set $K:=\{\sum a_{ij}x^iy^j\vert~ a_{ij}\in P(b)\}$. Since $[x, y]=b\in Z(D), K$ is closed under the multiplicative operation in $D$, so, clearly it is a subring of $D$. Moreover, since $x$ and $y$ both are algebraic over a finite field $P(b)$, it is easy to see that $K$ is finite. Therefore, $K$ is a finite division ring and by Wedderburn's Theorem  $K$ is a field. In particular, $x$ and $y$ commute with each other and that is a contradiction.\dpcm\\

\newpage
\noindent
{\em Proof of Theorem 1.}  Clearly we can suppose that $D$ is non-commutative. In the first, we shall prove that $G$ is abelian. Thus, suppose that $G$ is non-abelian. Then, there exist elements $x, y\in G$ such that $[x, y]\neq 1$. Let $H=\langle x, y\rangle $  be the subgroup of $G$ generated by $x$ and $y$. Set
$$H_0=H, H_1=[H_0, H_0]\mbox{ and } H_i=[H_{i-1}, H]\mbox{ for all } i\geq 2.$$

Since $H$ is non-abelian nilpotent subgroup, there exists the integer  $s\geq 1$ such that $H_s\neq 1$ and $H_{s+1}=1$. By definition we have
$$H_s=\langle[a, b]\vert~ a\in H_{s-1}, b\in H\rangle\neq 1.$$ 

Therefore, there exist $a_0\in H_{s-1}, b_0\in H$ such that $c:=[a_0, b_0]\neq 1$. Since $c\in H_s$ and $H_{s+1}=1, c$ commutes with each element from $H$. In particular, $c$ commutes with both $x$ and $y$. Now, suppose that 
$$G=G_n\triangleleft G_{n-1}\triangleleft\ldots\triangleleft G_1\triangleleft G_0=D^*.$$

By setting $D_1=C_D(c), N_i=G_i\cap D_1$ for all $i\in\{1, \ldots, n\}$, we obtain                                
$$N=N_n\triangleleft N_{n-1}\triangleleft\ldots\triangleleft N_1\triangleleft N_0=D_1^*.$$

So $N$ is a subnormal locally nilpotent subgroup of $D_1^*$. Since $x$ and $y$ are both commute with $c, H=\langle x, y\rangle\leq D_1^*$. It follows that $H\leq D_1^*\cap G_n=N_n=N$. Hence $a_0, b_0\in N$. Since $c=[a_0, b_0]\neq 1$, by Lemma 2, $c\not\in Z(D_1)$. In the other hand, since $D_1=C_D(c), c\in Z(D_1)$ and we have a contradiction. Thus, $G$ is abelian. Now, in view of Scott (1987, 14.4.4, p.440) $G$ is contained in $Z(D)$.\dpcm\\

\noindent
{\bf Acknowledgement.} The authors thank Mojtaba Ramezannassab for some his remarks concerning the proof of our theorem in Hai-Thin (2009, Th.2.2).

\end{document}